\documentclass[11pt, reqno]{amsart}
\usepackage{enumitem}
\usepackage{amsthm,amsmath,amssymb}

\usepackage{mathrsfs}

\numberwithin{equation}{section}

\newtheorem{theorem}{Theorem}[section]

\theoremstyle{plain}

\newtheorem{lemma}[theorem]{Lemma}

\newtheorem{remark}[theorem]{Remark}

\def\be{\begin{equation}}
	\def\ee{\end{equation}}

\def\d{\nabla}

\DeclareMathOperator{\A}{A}

\DeclareMathOperator{\Ric}{Ric}

\DeclareMathOperator{\vol}{Vol}

\DeclareMathOperator{\Hess}{Hess}

\numberwithin{equation}{section}

\begin{document}
	\title[]{A direct approach to sharp Li-Yau estimates on closed manifolds with negative Ricci lower bound}
	\author{ XingYu Song, Ling Wu and Meng Zhu}
	\address{School of Mathematical Sciences,  Key Laboratory of MEA (Ministry of Education) \& Shanghai Key Laboratory of PMMP,  East China Normal University, Shanghai 200241, China}
	\email{ 52215500013@stu.ecnu.edu.cn, 52215500012@stu.ecnu.edu.cn, mzhu@math.ecnu.edu.cn}
	\date{}

	\begin{abstract}
		Recently, Qi S.Zhang \cite{ZQ} has derived a sharp Li-Yau estimate for positive solutions of the heat equation on closed Riemannian manifolds with the Ricci curvature bounded below by a negative constant. The proof is based on an integral iteration argument which utilizes Hamilton's gradient estimate, heat kernel Gaussian bounds and parabolic Harnack inequality. 
		
		In this paper, we show that the sharp Li-Yau estimate can actually be obtained directly following the classical maximum principle argument, which simplifies the proof in \cite{ZQ}. In addition, we apply the same idea to the heat and conjugate heat equations under the Ricci flow and prove some Li-Yau type estimates with optimal coefficients.
		
	\end{abstract}
	\maketitle
	%\tableofcontents
	
	\section{Introduction}
	The renowned Li-Yau estimate \cite{LY} states that on an $n$-dimensional complete Riemannian manifold $(M^n,g)$ with the Ricci curvature $\Ric\geq -K$, $K\geq0$, any positive solution $u$ of the heat equation
	\begin{equation}\label{equ1.1}
		\left(\Delta-\frac{\partial}{\partial t}\right)u=0
	\end{equation}
	satisfies
	\be\label{equ1.2}
	\frac{|\nabla u|^2}{u^2}- \alpha \frac{\partial_t u}{u}\le \frac{n\alpha ^2}{2t}+\frac{n\alpha ^2K}{2(\alpha-1)}, \ \forall \alpha>1, \ t>0.
	\ee
	Although \eqref{equ1.2} was discovered earlier by Aronson-B\'enilan on Euclidean spaces \cite{AB}, Li and Yau were the first to prove it on general differentiable manifolds and demonstrated numerous crucial applications in the research of differential geometry and geometric analysis. The Li-Yau estimate is closely related to parabolic Harnack inequality, Sobolev and Poincar\'e inequalities, volume doubling property, heat kernel estimates, estimates of the Green's functions, eigenvalue estimates, etc. Even more powerfully, when the Ricci curvature is nonnegative, i.e., $K=0$, one can let $\alpha$ in \eqref{equ1.2} tend to $1$ and obtain the sharp Li-Yau estimate
	\be\label{equ1.3}
	\frac{|\nabla u|^2}{u^2}- \frac{\partial_t u}{u}\le \frac{n}{2t},
	\ee
	which in term gives an optimal heat kernel estimate, as well as the laplacian and Bishop-Gromov volume comparison theorems.
	
	Due to its significancy, many efforts have been made to improve or generalize Li-Yau estimates. Assuming $\Ric$ bounded below by a negative constant, the coefficient $\alpha$ in \eqref{equ1.2}  was improved in e.g. \cite{Yau, Ha, LX, Qian, BBG}.  The Li-Yau estimate was also generalized to various linear and semilinear parabolic equations, as well as matrix versions (see e.g. \cite{Ha, CN, Lij, CFL}). Li-Yau type estimates under time-dependent metrics, for example the Ricci flow and K\"ahler-Ricci flow, were studied in the past two decades, and are very useful tools in characterizing the singularities of the Ricci flow (see e.g. \cite{Ha2, Ha4, CHD, BCP, KZ, Per, CxZz, CH, LZ}). Furthermore, the Ricci lower bound assumption for Li-Yau estimate has recently been relaxed to integral Ricci conditions and Bakry-\'Emery Ricci curvature condition, or even the $RCD^*(K, N)$ spaces (see e.g. \cite{ZZ1, ZZ2, BL, GM, QZZ, ZhZx, Liy}). 
	
	Estimate \eqref{equ1.3} is sharp in the sense that the coefficient of $\frac{\partial_t u }{u }$ is equal to $1$, and the equality is achieved by the heat kernel $(4\pi t)^{-n/2}exp\{-\frac{|x-y|^2}{4t}\}$ on $\mathbb{R}^n$. However, estimate \eqref{equ1.2} is not sharp for $K>0$. For example, Hamilton \cite{Ha} showed that 
	\begin{equation}\label{a}
		\frac{|\nabla u|^2}{u^2}-e^{2Kt}\frac{\partial_t u}{u}\le e^{4Kt}\frac{n}{2t},
	\end{equation} 
	where the constant $\alpha>1$ is improved to an exponential function approaching $1$ as $t$ tending to $0$.
	Indeed, an open question asks if one can take $\alpha=1$ even for $K>0$. 
	
	Very recently, Qi S. Zhang \cite{ZQ} has invented a method which answers the above question affirmatively for closed manifolds. More precisely, he showed that 
	
	\begin{theorem}{(\cite{ZQ})}\label{th1.1}
		Let $(M^n,g)$ be an $n$-dimensional closed Riemannian manifold with $\Ric\ge -K $ for some constant $K\ge0$, and $u$ a positive solution of the heat equation \eqref{equ1.1} on $M\times[0,+\infty)$. Then there exist constants $c_1$ and $c_2$ depending only on $n$ such that 
		\begin{equation}\label{sharp Li-Yau}
			t\left(\frac{|\nabla u|^2}{u^2}-\frac{\partial_t u}{u}\right) \le \frac{n}{2}+\sqrt{2nK(1+Kt)(1+t)}diam_M+\sqrt{K(1+Kt)(c_1+c_2K)t}
		\end{equation}
		for all $(x,t)\in M\times [0,+\infty)$, where $diam_M$ is the diameter of $M$.
	\end{theorem}
	
	The proof in \cite{ZQ} is mainly based on an integral iteration argument with dedicated estimates, in which Hamilton's gradient estimate for the heat equation \cite{Ha} , Gaussian bounds of the heat kernel and differential Harnack inequality are used. 
	
	In this paper, we observe that Theorem \ref{th1.1} can actually be obtained by directly following the classical maximum principle arguments, and the proof can be simplified. Our proof is inspired by the usage of Hamilton's gradient estimate and the heat kernel bounds in \cite{ZQ}.

	In Theorem \ref{th1.1}, the Laplace operator is associated with a fixed metric. We can also consider the heat equation and conjugate heat equation under the Ricci flow, in which the Laplace operator evolves as the metric changes in time. For the heat equation, we show that
	
	\begin{theorem}\label{th1.3}
		Let $M^n$ be an n-dimensional closed Riemannian manifold, and $g(x,t)$, $(x,t)\in M\times [0,T]$, a solution to the Ricci flow  $\frac{\partial}{\partial t}g(x,t)=-2\Ric(x,t)$ on $M$. Assume that for some constants $K_1,K_2\ge0$, $-K_1g(x,t)\le \Ric(x,t)\le K_2g(x,t)$, $\forall (x,t)\in M\times [0,T]$. Let $u$ be a smooth positive solution of the heat equation $\left(\Delta_t-\frac{\partial}{\partial t}\right)u=0$ on $M\times[0,T]$, where $\Delta_t$ is the Laplace operator of the metric $g(x, t)$. Then there exist constants $c_5$ and $c_6$ depending only on $n$ such that 
		\begin{equation}
			t\left(\frac{|\nabla u|^2}{u^2}-\frac{\partial_t u}{u}\right) \le \frac{n}{2}+nK_3t +\sqrt{c_5K_1[e^{c_6K_3t}(t+d^2_0(M))]}\nonumber
		\end{equation}
		for all $(x,t)\in M\times [0,T]$, where $d_t(M)$ is the diameter of $(M^n, g(x,t))$ and $K_3=\max\{K_1,K_2\}$.
	\end{theorem}
	\begin{remark}
		When $K_1=0$, Theorem \ref{th1.3} reduces to Theorem 2.9 in \cite{BCP}.
	\end{remark}

	In Theorem \ref{th1.3}, the Ricci curvature is assumed to be uniformly bounded up to time $T$ before the singular time. The main reason is that the heat kernel Gaussian bounds may not hold uniformly till the singular time. This assumption can be relaxed to a so-called "Type-I" Ricci curvature bound, if we assume in addition that Perelman's $\nu$ entropy is bounded from below. 
	
	Recall that Perelman's $\mathcal{W}$ functional is defined as
	\begin{equation}
		\mathcal{W}(g,f,\tau)=(4\pi\tau)^{-\frac{n}{2}}\int_M(\tau(R+|\nabla f|^2)+f-n)e^{-f}dg, \nonumber
	\end{equation}
	where $f\in C^{\infty}(M)$ and $\tau>0$. Then the $\mu$ and $\nu$ entropy functionals are 
	\begin{equation}
		\mu [g,\tau]:=\inf \left\{\mathcal{W}(g,f,\tau);\ f\in C^{\infty}(M)\ and\ \int_M(4\pi\tau)^{-\frac{n}{2}}e^{-f}dg=1\right\},\nonumber
	\end{equation}
	and
	\begin{equation}
		\nu[g,\tau]:=\inf\limits_{s\in[0,\tau]}\mu [g,s] .\nonumber
	\end{equation}
	
	\begin{theorem}\label{th1.5}
		Let $M^n$ be an n-dimensional closed Riemannian manifold and $g(x,t), (x,t)\in M\times [0,2)$, a solution to the Ricci flow  $\frac{\partial}{\partial t}g(x,t)=-2\Ric(x,t)$ on $M$. Assume that for some constant $K\ge0$, $-\frac{K}{2-t}g(x,t)\le \Ric(x,t)\le \frac{K}{2-t}g(x,t)$ and $\nu[g_0,4]\ge -nK $, $\forall (x,t)\in M\times [0,2)$. Let $u$ be a smooth positive solution of  the heat equation $\left(\Delta_t-\frac{\partial}{\partial t}\right)u(x,t)=0$ on $M\times [0,2)$. Then there exists a constant $c_9$ depending on $n,K$ such that 
		\begin{equation}
			t\left(\frac{|\nabla u|^2}{u^2}-\frac{\partial_t u}{u}\right) \le \frac{n}{2}+ \frac{nKt}{2-t}+\sqrt{\frac{Kc_9}{2-t}\left(t+d_0^2(M)t+d_0^2(M)\right)}\nonumber
		\end{equation}
		for all $(x,t)\in M\times [0,2)$, where $d_0(M)$ is the diameter of $(M,g(x,0))$.
	\end{theorem}
	
	Finally, we consider the Li-Yau estimate for positive solutions of the conjugate heat equation under the Ricci flow. 
	
	\begin{theorem}\label{th1.9}
		Let $M^n$ be an n-dimensional closed Riemannian manifold and $g(x,t), (x,t)\in M\times [0,T)$, a solution to the Ricci flow  $\frac{\partial}{\partial t}g(x,t)=-2\Ric(x,t)$ on $M$. Assume that for some constants $K_1,K_2\ge0$, $-K_1g(x,t)\le \Ric(x,t)\le K_2g(x,t)$, $|\nabla R|(x,t)\le K_2^{\frac{3}{2}}$ and $R(x,t)\ge0$,  $\forall (x,t)\in M\times [0,T]$. Let $u$ be a positive solution of the conjugate heat equation $\left(\Delta_t+\frac{\partial}{\partial t}-R\right)u=0$ on  $M\times [0,T]$. Then there exist constants $c_{10}$ and $ c_{11}$ depending only on $n$ such that 
		\begin{equation}
			\frac{|\nabla u|^2}{u^2}+\frac{u_{t}}{u}\le \frac{n}{T-t} +c_{10}e^{c_{11}K_1T}\left(1+\frac{d_0^2(M)}{T}\right)^2\left(\frac{1}{t}+K_1+K_2\right)\nonumber
		\end{equation}
		for all $(x,t)\in M\times (0, \frac{T}{2}]$, where $d_0(M)$ is the diameter of $(M,g(x,0))$.
	\end{theorem}

	The rest of the paper is organized as follows. In section 2, we give an alternative proof of Theorem \ref{th1.1}. Section 3 is devoted to the proof of Theorem \ref{th1.3}. Theorem \ref{th1.5} is dealt with separately in section 4, since the curvature is blowing up at the singular time and a different heat kernel estimate is required. Finally, we show Theorem \ref{th1.9} in section 5.

	\section{The heat equation with fixed metric}
	In this section, we prove the sharp Li-Yau estimate \eqref{sharp Li-Yau} .\\
	
	\noindent{\it Proof of Theorem \ref{th1.1}.} Let $u=u(x,t)$ be a positive solution of the heat equation on $M\times[0,+\infty)$. Denote 
	\begin{equation}\label{s1}
		Q=\frac{|\nabla u|^2}{u^2}-\frac{\partial_t u}{u}=|\d \ln u|^2-\partial_t\ln u=-\Delta \ln u. 
	\end{equation}
	%By direct computation, we have 
	%\begin{equation}
	%	\Delta \ln u=\frac{\partial_t u}{u}-\frac{|\nabla u|^2}{u^2}.\nonumber
	%\end{equation}
	From the Bochner formula and the Cauchy-Schwarz inequality, it follows that
	\begin{equation}\label{eq1.1}
		\begin{aligned}
			(\Delta-\partial_t)Q+2\left<\nabla\ln u,\nabla Q\right>&=2\left[|\Hess \ln u|^2+\Ric(\nabla \ln u,\nabla \ln u)\right]\\
			%&\ge 2 \left[\frac{(\Delta \ln u)^2}{n}-K|\nabla \ln u|^2\right]\\
			&\geq 2\left[\frac{Q^2}{n}-K|\nabla \ln u|^2\right].
			%&=2\left(\frac{1}{n}Q^2-\frac{2(\alpha-1)}{n}|\nabla \ln u|^2Q+\frac{(\alpha-1)^2|\nabla \ln u|^4}{n}-K|\nabla\ln u |^2\right),
		\end{aligned}
	\end{equation}
	Hence, 
	\begin{equation}\label{eq1.2}
		(\Delta-\partial_t)(tQ)+2\left<\nabla\ln u,\nabla (tQ)\right>\ge 2t\left(\frac{1}{n}Q^2-K|\nabla\ln u |^2\right)-Q.
	\end{equation}
	For any $T>0$, the maximum of $tQ$ on $M\times[0,T]$ is achieved at some $(x_0, t_0)$. We may assume that $tQ$ is positive at $(x_0,t_0)$, for otherwise the theorem is trivial. Thus, $t_0>0$, and at $(x_0,t_0)$ we have
	\begin{equation}
		\nabla(tQ)=0,\ \frac{\partial}{\partial t}(tQ)\ge0,\ and\  \Delta(tQ)\le0.\nonumber
	\end{equation}
	Evaluating \eqref{eq1.2} at $(x_0,t_0)$, we obtain 
	\begin{equation}
		0\ge\frac{2t_0}{n}Q^2(x_0,t_0)-Q(x_0,t_0)-2Kt_0|\nabla\ln u|^2(x_0,t_0).\nonumber
	\end{equation}
	Multiplying both sides above by $t_0$ gives
	\begin{equation}
		0\ge \frac{2}{n}(t_0Q(x_0,t_0))^2-t_0Q(x_0,t_0)-2Kt_0^2|\nabla\ln u|^2(x_0,t_0),\nonumber
	\end{equation}
	which infers 
	\begin{equation}\label{eq1.3}
		t_0Q(x_0,t_0)\le \frac{n}{2}+\sqrt{nK}t_0|\nabla\ln u|(x_0,t_0).
	\end{equation}
	
	Next we will bound $t_0|\nabla\ln u|(x_0,t_0)$. From \cite{ZQ}, the heat kernel $H(x,t,y)$, $x,y\in M$, $t>0$ satisfies the estimates
	\begin{equation}\label{eq1.4}
		t|\nabla_x\ln H(x,t,y)|^2\le 2(1+Kt)\left(2C_1Kt+4C_2^2K+\frac{diam_M^2}{t}+2\ln C_3\right),\ \forall\  0<t\le 8,
	\end{equation}
	and
	\begin{equation}\label{eq1.5}
		t|\nabla_x\ln H(x,t,y)|^2\le2(1+Kt)(n\ln2+2C_4K+diam_M^2),\ \forall \ t\ge 8,
	\end{equation}
	where the positive constants $C_i$, $1\le i\le 4$, depend only on $n$. Actually, \eqref{eq1.4} comes from Hamilton estimate \cite{Ha} and Gaussian bounds of the heat kernel \cite{LY}, and \eqref{eq1.5} relies further on the parabolic Harnack inequality.
	
	Combining \eqref{eq1.4} and \eqref{eq1.5}, we deduce that
	\begin{equation}\label{1.6}
		t^2|\nabla_x\ln H(x,t,y)|^2\le 2(1+KT)(C_5T+C_6KT+diam_M^2+diam_M^2T),\ \forall\ 0<t\le T.
	\end{equation}
	Since $T$ is arbitrary, we conclude that
	\begin{equation}\label{eq1.7}
		t^2|\nabla_x\ln H(x,t,y)|^2\le 2(1+Kt)(C_5t+C_6Kt+diam_M^2+diam_M^2t),\ \forall\ t>0\ and\ x,y\in M.
	\end{equation}
	A short argument from \cite{YZ} implies that the bound \eqref{eq1.7} also holds if one replaces the heat kernel by any positive solution $u(x,t)$ of the heat equation. For completeness, we present the proof below. By Duhamel's principle, we have $u(x,t)=\int_M H(x,t,y)u(y,0)dy$. Then
	\begin{equation}
		\begin{aligned}
			|\nabla_x u(x,t)|^2&=\int_M\int_M \left<\nabla_xH(x,t,y),\nabla_xH(x,t,z)\right>u(y,0)u(z,0)dydz\\
			&\le\frac{1}{2}\int_M\int_M\left(\frac{|\nabla_xH(x,t,y)|^2}{H(x,t,y)}H(x,t,z)+\frac{|\nabla_xH(x,t,z)|^2}{H(x,t,z)}H(x,t,y)\right)u(y,0)u(z,0)dydz.\nonumber
		\end{aligned}
	\end{equation}
	Using \eqref{eq1.7} in the integrand, one has 
	\begin{equation}\label{eq1.8}
		t^2|\nabla \ln u(x,t)|^2\le 2(1+Kt)(C_5t+C_6Kt+diam_M^2+diam_M^2t),\ \forall\ t\ge0\ and\ x\in M.
	\end{equation}
	Plugging \eqref{eq1.8} into \eqref{eq1.3}, we obtain 
	\begin{equation}\label{eq1.9}
		TQ(x,T)\le t_0Q(x_0,t_0)\le \frac{n}{2}+\sqrt{K(1+KT)(C_7T+C_8KT)}+\sqrt{2nK(1+KT)(1+T)}diam_M.
	\end{equation}
	By the arbitrariness of $T$ again, Theorem \ref{th1.1} follows from \eqref{eq1.9}.  \qed 
	
	\begin{remark}
	One may replace \eqref{eq1.3} by the Li-Yau type inequalities of Yau \cite{Yau} or Bakry-Qian \cite{BQ} and obtain a similar sharp Li-Yau estimate. However, in the proof above we only need this type of inequality at a maximal point, and \eqref{eq1.3} is easier to derive. 
	\end{remark}
	
	\begin{remark}
		If $\Ric\ge0$, then Theorem \ref{th1.1} becomes the sharp Li-Yau estimate \eqref{equ1.3}.
	\end{remark}

	\begin{remark}\label{rm1.2}
		Under the same condition in Theorem \ref{th1.1}, for any $\alpha>1$, using Theorem \ref{th1.1} and \eqref{eq1.8}, one gets that there exist constants $c_i$ $(i=1,2,3,4)$ depending only on $n$ such that
		\begin{equation}
			\begin{aligned}
				t\left(\alpha\frac{|\nabla u|^2}{u^2}-\frac{u_t}{u}\right)&=t\left(\frac{|\nabla u|^2}{u^2}-\frac{u_t}{u}\right)+(\alpha-1)t\frac{|\nabla u|^2}{u^2}\\
				&\le \frac{n}{2}+\sqrt{2nK(1+Kt)(1+t)}diam_M+\sqrt{K(1+Kt)(c_1+c_2K)t}\\
				&\ \ \ \ +2(\alpha-1)\frac{1+Kt}{t}\left(c_3t+c_4Kt+diam_M^2+diam_M^2t\right)\nonumber
			\end{aligned}
		\end{equation}
		for all $(x,t)\in M\times (0,+\infty)$.
	\end{remark}

	\section{The heat equation under the Ricci flow with uniform Ricci bound}
	
	In this section, we consider Li-Yau estimates in the Ricci flow case and prove Theorem \ref{th1.3}. First, we need the following Hamilton type estimate by Qi S.Zhang \cite{ZQ2} and X. Cao-R. Hamilton \cite{CH}.
	
	\begin{lemma}\label{lm3.1}(\cite{ZQ2}, \cite{CH})
		Suppose that $M^n$ is an $n$-dimensional closed manifold, and $g(x,t)$, $t\in [0,T]$ is a solution to the Ricci flow  $\frac{\partial}{\partial t}g(x,t)=-2\Ric(x,t)$ on $M$. Let $u$ be a positive solution of the heat equation $\left(\Delta_{t}-\frac{\partial}{\partial t}\right)u(x,t)=0$ on $ M\times [0,T]$. Then it holds that
		\begin{equation}
			\frac{|\nabla u|}{u}\le \sqrt{\frac{1}{t}\ln \frac{\A}{u}}\nonumber
		\end{equation} 
		for all $(x,t)\in M\times (0,T] $ with $\A=\sup\limits_{x\in M}u(x,0)$.
	\end{lemma}
	In \cite{ZM}, the third author provides Gaussian upper and lower bounds for the heat kernel of the heat equation under the Ricci flow.
	\begin{lemma}\label{lm3.2}(\cite{ZM})
		Let $(M^n,g(x,t))_{t\in [0,T)}$ be a complete solution to the Ricci flow  $\frac{\partial}{\partial t}g(x,t)=-2\Ric(x,t)$, $(x,t)\in M\times [0,T)$ for $T < \infty$. Suppose that $\Ric(x,t)\ge-K_1g(x,t)$ on $M\times [0,T)$ for some nonnegative constant $K_1$ and $\Lambda=\int_{0}^{T}\sup\limits_{x\in M}|\Ric|(x,t)dt< \infty$. Let $H(x,t;y,l)$ be the heat kernel of the heat equation $\left(\Delta_{t}-\frac{\partial}{\partial t}\right)u(x,t)=0$ with $0\le l<t<T$. Then we have the following upper and lower bounds:
		\begin{equation}
			H(x,t;y,l)\le \tilde{c}_1(n)e^{\tilde{c}_2(n)e^{{\tilde{c}_3(n)\Lambda+{\tilde{c}_4(n)K_1T}}}}\min\left\{\frac{e^{-\frac{d_t^2(x,y)}{8e^{4K_1T}(t-l)}}}{\vol_l(B_l(x,\sqrt{\frac{t-l}{8}}))},\ \frac{e^{-\frac{d_t^2(x,y)}{8e^{4K_1T}(t-l)}}}{\vol_t(B_t(y,\sqrt{\frac{t-l}{8}}))} \right\}\nonumber
		\end{equation}
		and
		\begin{equation}
			H(x,t;y,l)\ge \frac{\tilde{c}_5(n)e^{-\tilde{c}_6(n)e^{\tilde{c}_7(n)\Lambda+\tilde{c}_8(n)K_1T}}e^{-\frac{4d_t^2(x,y)}{t-l}}}{\vol_t(B_t(y,\sqrt{\frac{t-l}{8}}))},\nonumber
		\end{equation}
		where $B_t(x,\sqrt{\frac{t-l}{8}})$ is the geodesic ball of radius $\sqrt{\frac{t-l}{8}}$ in $M$ at time $t$ and, $\vol_t(.)$ denotes the volume of a region at time $t$.
	\end{lemma}
	\noindent{\it Proof of Theorem \ref{th1.3}.} Denote by $f=\ln u$ and $F=t(|\nabla f|^2-f_t)=-t\Delta f$. By direct computation, we obtain 
	\begin{equation}
		\Delta F =t(2|\Hess f|^2+2\left<\nabla f,\nabla \Delta f\right>+2\Ric(\nabla f,\nabla f)-(\Delta f)_t+2\left<\Ric,\Hess f\right>),\nonumber
	\end{equation}
	where we have used the identity 
	\begin{equation}
		\Delta(f_t)=(\Delta f)_t-2\left<\Ric,\Hess f\right>.\nonumber
	\end{equation}
	On the other hand, 
	\begin{equation}
		F_t=(|\nabla f|^2-f_t)+t(2\left<\nabla f,\nabla f_t\right>+2\Ric(\nabla f,\nabla f)-f_{tt}).\nonumber
	\end{equation}
	%where we noticed 
	%\begin{equation}
	%	\partial_t|\nabla f|^2=2\left<\nabla f,\nabla f_t\right>+2\Ric(\nabla f,\nabla f).\nonumber
	%\end{equation}
	Since  $(\Delta -\partial_t) f=-|\nabla f|^2$, we have
	\begin{equation}
		\begin{aligned}
			2t&\left<\nabla f,\nabla \Delta f\right>-t(\Delta f)_t-2 t\left<\nabla f,\nabla f_t\right>+tf_{tt}\\
			&=2 t\left<\nabla f,\nabla(\Delta f-f_t)\right>+t(f_t-\Delta f)_t\\
			&=-2 t\left<\nabla f,\nabla |\nabla f|^2\right>+t(|\nabla f|^2)_t\\
			&=-2 \left<\nabla f,\nabla(F+tf_t)\right>+2t\left<\nabla f,\nabla f_t\right>+2t\Ric(\nabla f,\nabla f)\\
			&=-2\left<\nabla f,\nabla F\right>+2t\Ric(\nabla f,\nabla f). \nonumber
		\end{aligned}
	\end{equation}
	Therefore,
	\begin{equation}\label{e3.1}
		(\Delta-\partial_t)F=-2\left<\nabla f,\nabla F\right>+2t( |\Hess f|^2+\left<\Ric, \Hess f \right>)+2t\Ric(\nabla f, \nabla f)-(|\nabla f|^2-f_t).
	\end{equation}
	For any $0<a<1$, by the Cauchy-Schwarz inequality and the bound of the Ricci curvature, one has
	\begin{equation}\label{e3.2}
		\begin{aligned}
			|\Hess f|^2+\left<\Ric, \Hess f \right>&=\sum_{i,j=1}^{n}( f^2_{ij}+R_{ij}f_{ij})\\
			&\ge \sum_{i,j=1}^{n}\left( f_{ij}^2-(1-a)f^2_{ij}-\frac{1}{4(1-a)}R_{ij}^2\right)\\
			&\ge\sum_{i,j=1}^{n}af_{ij}^2-\frac{n}{4(1-a)}K_3^2\\
			%&=a|\Hess f|^2-\frac{n}{4(\alpha-a)}K_3^2\\
			&\ge \frac{a}{n}(\Delta f)^2-\frac{n}{4(1-a)}K_3^2\\
			&=\frac{a}{n}(f_t-|\nabla f|^2)^2-\frac{n}{4(1-a)}K_3^2.
		\end{aligned}
	\end{equation}
	Plugging \eqref{e3.2} into \eqref{e3.1}, we see that 
	\begin{equation}\label{e3.4}
		(\Delta-\partial_t)F+2\left<\nabla f,\nabla F\right>\ge \frac{2a}{n}\frac{F^2}{t}-\frac{F}{t}-\frac{nt}{2(1-a)}K_3^2-2K_1t|\nabla f|^2
	\end{equation}
	for any $0<a<1$.\\
	For any $0<\tau\le T$, let $(x_0,t_0)\in M\times[0,\tau]$ be a maximum point of $F$. Without loss of generality, we may assume $F(x_0,t_0)>0$. Then $t_0>0$, and by the maximum principle, we have 
	\begin{equation}\label{e3.5}
		\Delta F(x_0,t_0)\le0,\ \frac{\partial}{\partial t}F(x_0,t_0)\ge0,\ and\ \nabla F(x_0,t_0)=0.
	\end{equation}
	It follows from \eqref{e3.4} and \eqref{e3.5} that at $(x_0,t_0)$,
	\begin{equation}
		\frac{2a}{n}F^2(x_0,t_0)-F(x_0,t_0)-\frac{nt_0^2}{2(1-a)}K_3^2-2K_1t_0^2|\nabla f|^2(x_0,t_0)\le0,\nonumber
	\end{equation}
	which yields
	\begin{equation}
		F(x_0,t_0)\le \frac{n}{4a}\left(1+\sqrt{1+\frac{4at_0^2}{1-a}K_3^2}\right)+\sqrt{\frac{n}{a}K_1}t_0|\nabla f|(x_0,t_0).\nonumber
	\end{equation}
	By choosing $a=\frac{1+K_3t_0}{1+2K_3t_0}\in [\frac{1}{2},1)$, one has
	\begin{equation}\label{e3.6}
		F(x_0,t_0)\le nK_3t_0+\frac{n}{2}+\sqrt{2nK_1}t_0|\nabla f|(x_0,t_0).
	\end{equation}
	Again, to estimate $t|\d f|$, we first find a bound for $t|\nabla \ln H(x,t;y,0)|$ with $x,y\in M$ and $t\in(0,\tau]$, where $H(x,t;y,0)$ is the heat kernel of the heat equation $\left(\Delta_{t}-\frac{\partial}{\partial t}\right)H(x,t;y,0)=0$.
	
	For a time $0<t_1\le\frac{\tau}{2}$, we denote by $h(x,t)=H(x,t+t_1;y,0)$ with $x,y\in M$ and $t\in[0,t_1]$.\\
	Since $h(x,t)$ is a positive solution of the heat equation, by Lemma \ref{lm3.1}, we derive 
	\begin{equation}\label{e3.7}
		t|\nabla \ln H(x,t+t_1;y,0)|^2\le \ln \frac{A}{H(x,t+t_1;y,0)},\ (x,t)\in M\times(0,t_1],
	\end{equation}
	where $A=\sup\limits_{x\in M}H(x,t_1;y,0)$.
	
	According to Lemma \ref{lm3.2}, since $\Lambda\le 2\sqrt{n}K_3t_1$, the following upper and lower bounds for $h(x,t)$ hold
	\begin{equation}
		\frac{\tilde{c}_5(n)e^{-\tilde{c}_6(n)e^{\tilde{c}_{10}(n)K_3t_1}}e^{-\frac{4d^2_{t+t_1}(M)}{t_1}}}{\vol_{t+t_1}(B_{t+t_1}(y,\sqrt{\frac{t_1}{4}}))}		\le h(x,t) \le \frac{\tilde{c}_1(n)e^{\tilde{c}_2(n)e^{\tilde{c}_9(n)K_3t_1}}}{\vol_{t+t_1}(B_{t+t_1}(y,\sqrt{\frac{t+t_1}{8}}))},\ t\in[0,t_1].\nonumber
	\end{equation}
	The upper bound gives 
	\begin{equation}
		A=\sup\limits_{x\in M}H(x,t_1;y,0) \le \frac{\tilde{c}_1(n)e^{\tilde{c}_2(n)e^{\tilde{c}_9(n)K_3t_1}}}{\vol_{t_1}(B_{t_1}(y,\sqrt{\frac{t_1}{8}}))},\nonumber
	\end{equation}
	which, combining with the lower bound of $h(x,t)$, yields that
	\begin{equation}\label{e3.8}
		\frac{A}{H(x,t+t_1;y,0)}\le \tilde{c}_{11}(n)e^{\tilde{c}_{12}(n)e^{\tilde{c}_{13}(n)K_3t_1}}e^{\frac{4d^2_{t+t_1}(M)}{t_1}} \frac{\vol_{t+t_1}(B_{t+t_1}(y,\sqrt{\frac{t_1}{4}}))}{\vol_{t_1}(B_{t_1}(y,\sqrt{\frac{t_1}{8}}))}.
	\end{equation}
	
	The uniform boundedness of the Ricci curvature implies that $g(x,t)$ is uniformly equivalent to the initial metric $g(0)$ (see Corollary 6.11 in \cite{BLN} ), that is ,
	\begin{equation}\label{e3.9}
		e^{-2K_2t}g(0)\le g(x,t) \le e^{2K_1t}g(0).
	\end{equation}
	Then, according to the Bishop-Gromov volume comparison theorem, we have 
	\begin{equation}
		\begin{aligned}
			\frac{\vol_{t+t_1}(B_{t+t_1}(y,\sqrt{\frac{t_1}{4}}))}{\vol_{t_1}(B_{t_1}(y,\sqrt{\frac{t_1}{8}}))}&\le\frac{e^{3nt_1K_3}\vol_{t_1}(B_{t_1}(y,\sqrt{\frac{t_1}{4}}))}{\vol_{t_1}(B_{t_1}(y,\sqrt{\frac{t_1}{8}}))}\\
			&\le e^{3nt_1K_3}(\sqrt{2})^ne^{\sqrt{(n-1)K_3t_1}}\\
			&\le e^{\tilde{c}_{14}(n)(K_3t_1+1)}.\nonumber
		\end{aligned}
	\end{equation}
	Plugging this into \eqref{e3.8}, we see that  \eqref{e3.7} becomes
	\begin{equation}
		t|\nabla \ln H(x,t+t_1;y,0)|^2 \le \tilde{c}_{15}(n)\left(e^{\tilde{c}_{16}(n)K_3t_1}+\frac{d^2_{t+t_1}(M)}{t_1}\right),\ \forall t\in (0,t_1].\nonumber
	\end{equation}
	Taking $t=t_1$ and using the arbitrariness of $t_1$, one concludes that
	%\begin{equation}
	%		t|\nabla \ln H(x,t;y,0)|^2\le\tilde{c}_{17}(n)\left(e^{\tilde{c}_{18}(n)K_3t}+\frac{d^2_{t}(M)}{t}\right), \forall\ %0<t\le\tau, \nonumber
	%\end{equation}
	%which infers
	\begin{equation}\label{e3.10}
		t^2|\nabla \ln H(x,t;y,0)|^2\le\tilde{c}_{17}(n)\left(e^{\tilde{c}_{19}(n)K_3t}(t+d^2_0(M))\right), \forall\ 0<t\le\tau.
	\end{equation}
	Since the solution $u$ to the heat equation can also be written as $$u(x,t)=\int_M H(x,t;y,0)u(y,0)d\mu_0(y),$$ where $d\mu_0$ is the volume element at $t=0$, the same bound as in \eqref{e3.10} actually holds for any positive solution of the heat equation, that is,
	\begin{equation}\label{e3.11}
		t^2|\nabla \ln u(x,t)|^2 \le \tilde{c}_{17}(n)\left(e^{\tilde{c}_{19}(n)K_3t}(t+d^2_0(M))\right), \forall\ 0\le t\le\tau.
	\end{equation}
	Combining \eqref{e3.11} with \eqref{e3.6}, we obtain 
	\begin{equation}\label{e3.12}
		F(x,\tau)\le F(x_0,t_0)\le nK_3\tau +\frac{n}{2}+\sqrt{\tilde{c}_{20}(n)K_1[e^{\tilde{c}_{19}(n)K_3\tau}(\tau+d^2_0(M))]}.
	\end{equation}
	Since $\tau\in (0,T]$ can be chosen arbitrarily, it finishes the proof of Theorem \ref{th1.3}.  \qed
	
	\begin{remark}\label{rm1.4}
		Under the assumptions and notations as in Theorem \ref{th1.3}, similar to Remark \ref{rm1.2}, we have for any $\alpha>1$, there exist constants $c_i$ $(i=5,6,7,8)$ depending only on $n$ such that 
		\begin{equation}
			\begin{aligned}
				t\left(\alpha\frac{|\nabla u|^2}{u^2}-\frac{\partial_t u}{u}\right)&\le \frac{n}{2}+nK_3t +\sqrt{c_5K_1[e^{c_6K_3t}(t+d^2_0(M))]}\\
				&\ \ \ \ +(\alpha-1)\frac{1}{t}c_7e^{c_8K_3t}(t+d_0^2(M))\nonumber
			\end{aligned}
		\end{equation}
		for all $(x,t)\in M\times (0,T]$.
	\end{remark}
	
	\section{The heat equation under the Ricci flow with blowing-up Ricci curvature}
	
	To prove Theorem \ref{th1.5}, we will need the following Gaussian upper and lower bounds for the heat kernel near the singular time.
	\begin{lemma}\label{lm4.1}(\cite{HM})
		Let $(M^n,g(x,t))_{t\in [-2,0)}$ be a closed solution of Ricci flow satisfying $\nu[g_{-2},4]\ge-K$ and $|R(x,t)|\le K|t|^{-1}$ or some constant $K\ge0$. Let $H(x,t;y,l)$ be the heat kernel of the equation $\left(\Delta_{t}-\frac{\partial}{\partial t}\right)H(x,t;y,l)=0$ with $-2\le l<t<0$ and $x,y\in M$. Then, for any $x,y\in M$ and $-\frac{1}{2}\le l <t<0$, there exists a constant $C^{*}=C^{*}(K,n)$ such that
		\begin{equation}
			\frac{1}{C^*(t-l)^{\frac{n}{2}}}e^{-\frac{C^*d^2_l(x,y)}{t-l}}\le H(x,t;y,l)\le\frac{C^*}{(t-l)^{\frac{n}{2}}}e^{-\frac{d^2_l(x,y)}{C^*(t-l)}}.\nonumber
		\end{equation}
	\end{lemma}
\medskip

	\noindent{\it Proof of Theorem \ref{th1.5}.} 
	Denote $f=\ln u$ and $F=t(|\nabla f|^2-f_t)$.
	
	For any $0<T<2$, let $(x_0,t_0)\in M \times(0,T]$ be a maximum point of $F$ and $F(x_0,t_0)>0$. Similar to  \eqref{e3.6}, we get 
	\begin{equation}\label{e4.1}
		F(x_0,t_0)\le nt_0\frac{K}{2-t_0}+\frac{n}{2} +\sqrt{\frac{2nK}{2-t_0}}t_0|\nabla f|(x_0,t_0)\ \  for\ \ 0< t_0\le T. 
	\end{equation}
	Since the Ricci curvature bound implies $|R(x,t)|\le\frac{nK}{2-t}$, combining the condition $\nu[g_0,4]\ge-nK$, by applying Lemma \ref{lm4.1}, there exists a constant $\bar{C_1}(n,K)$ such that for $\frac{3}{2}<t<2$ and any $x,y\in M$,
	\begin{equation}\label{e4.2}
		\frac{1}{\bar{C}_1(n,K)\left(t-\frac{3}{2}\right)^{\frac{n}{2}}}e^{-\frac{\bar{C_1}(n,K)d^2_{\frac{3}{2}}(x,y)}{t-\frac{3}{2}}}\le H(x,t;y,\frac{3}{2})\le\frac{\bar{C}_1(n,K)}{\left(t-\frac{3}{2}\right)^{\frac{n}{2}}}e^{-\frac{d^2_{\frac{3}{2}}(x,y)}{\bar{C}_1(n,K)\left(t-\frac{3}{2}\right)}}.
	\end{equation}
	For a time $\frac{3}{2}<t_1<\frac{8}{5}$, we denote $h(x,t)=H(x,t+t_1;y,\frac{3}{2})$, $t\in\left[0,\frac{1}{4}t_1\right]$, then $t+t_1\in \left[t_1,\frac{5}{4}t_1\right]\subset \left(\frac{3}{2},2\right)$. Hence, \eqref{e4.2} implies 
	\begin{equation}
		\bar{C}_2(n,K)e^{-\frac{\bar{C}_1(n,K)d^2_{\frac{3}{2}}(M)}{t_1-\frac{3}{2}}}\le H(x,t+t_1;y,\frac{3}{2})\le\frac{\bar{C}_1(n,K)}{\left(t+t_1-\frac{3}{2}\right)^{\frac{n}{2}}}\ \ for\ \ t\in\left[0,\frac{1}{4}t_1\right].\nonumber
	\end{equation}
	The upper bound gives 
	\begin{equation}
		A=\sup\limits_{x\in M}H(x,t_1;y,\frac{3}{2})\le \frac{\bar{C}_1(n,K)}{(t_1-\frac{3}{2})^{\frac{n}{2}}}.\nonumber
	\end{equation}
	By lemma \ref{lm3.1},
	\begin{equation}
		t|\nabla \ln H(x,t+t_1;y,\frac{3}{2})|^2\le \bar{C}_3(n,K) +\ln(t_1-\frac{3}{2})^{-\frac{n}{2}}+\frac{\bar{C}_1(n,K)d^2_{\frac{3}{2}}(M)}{t_1-\frac{3}{2}}.\nonumber
	\end{equation}
	Since
	\begin{equation}
		-\frac{K}{2-t}g(x,t)\le \Ric(x,t)\le \frac{K}{2-t}g(x,t),\nonumber
	\end{equation}
	then
	\begin{equation}
		\left(\frac{2}{2-t}\right)^{-2K}g(0)\le g(x,t) \le \left(\frac{2}{2-t}\right)^{2K}g(0),\nonumber
	\end{equation}
	therefore
	\begin{equation}
		t|\nabla \ln H(x,t+t_1;y,\frac{3}{2})|^2\le \bar{C}_3(n,K) +\ln(t_1-\frac{3}{2})^{-\frac{n}{2}}+\frac{\bar{C}_4(n,K)d^2_{0}(M)}{t_1-\frac{3}{2}}.\nonumber
	\end{equation}
	We take $t=\frac{1}{4}t_1$ and use $t_1\in \left(\frac{3}{2},\frac{8}{5}\right)$ to conclude that 
	\begin{equation}
		\frac{1}{5}t|\nabla \ln H(x,t;y,\frac{3}{2})|^2\le \bar{C}_3(n,K) +\ln\left(\frac{4}{5}t-\frac{3}{2}\right)^{-\frac{n}{2}}+\frac{\bar{C}_4(n,K)d^2_{0}(M)}{\frac{4}{5}t-\frac{3}{2}},\ \ \forall t\in \left(\frac{15}{8},2\right),\nonumber
	\end{equation}
	in particular,
	\begin{equation}
		t|\nabla \ln H(x,t;y,\frac{3}{2})|^2\le\bar{C}_5(n,K)(1+d^2_0(M)),\ \ \forall t\in\left(\frac{17}{9},2\right).\nonumber
	\end{equation}
	Since  $u(x,t)=\int_M H(x,t;y,\frac{3}{2})u(y,\frac{3}{2})d\mu_{\frac{3}{2}}(y)$, $t\ge\frac{3}{2}$, the above upper bound actually implies 
	\begin{equation}\label{e4.3}
		t|\nabla \ln u(x,t)|^2\le \bar{C}_5(n,K)(1+d^2_0(M)),\ \ \forall t\in\left(\frac{17}{9},2\right).
	\end{equation}
	For $t\in \left[0,\frac{17}{9}\right]$, since the Ricci curvature is uniformly bounded, from \eqref{e3.11}, we derive 
	\begin{equation}\label{e4.4}
		t^2|\nabla \ln u(x,t)|^2\le \bar{C}_6(n)\left(e^{\bar{C}_7(n)Kt}(t+d_0^2(M))\right),\ \ \forall t\in\left[0,\frac{17}{9}\right].
	\end{equation}
	Combining \eqref{e4.3} and \eqref{e4.4}, we conclude that
	\begin{equation}
		t^2|\nabla \ln u(x,t)|^2\le\bar{C}_8(n,K)(t+d^2_0(M)t+d^2_0(M)).\ \ \forall t\in [0,2).\nonumber
	\end{equation}
	Plugging above into \eqref{e4.1} gives
	\begin{equation}
		F(x_0,t_0)\le\frac{n}{2}+\frac{nKt_0}{2-t_0}+\sqrt{\frac{K\bar{C}_9(n,K)}{2-t_0}(t_0+d_0^2(M)t_0+d_0^2(M))},\nonumber
	\end{equation}
	which implies
	\begin{equation}
		F(x,T)\le F(x_0,t_0)\le \frac{n}{2}+\frac{nKT}{2-T}+\sqrt{\frac{K\bar{C}_9(n,K)}{2-T}(T+d_0^2(M)T+d_0^2(M))}.\nonumber
	\end{equation}
	The arbitrariness of $T$ gives the theorem. \\ \qed
	
	\section{The Conjugate heat equation under the Ricci flow}
	In this section, we deal with the conjugate heat equation under the Ricci flow and prove Theorem \ref{th1.9}.
	
	First, we recall the following Hamilton type estimate for the positive solution to the conjugate heat equation under the Ricci flow.
	\begin{lemma}\label{lm5.1}(\cite{BLN})
		Let $(M^n,g(t))_{t\in [0,T]}$ be a solution to the Ricci flow $\frac{\partial}{\partial t}g(x,t)=-2\Ric(x,t)$. For a fixed point $o\in M$, assume that for some constants $K_1,\ K_2\ge 0$,  $\Ric(x,t)\ge -K_1g(x,t)$ and $|\nabla R|(x,t)\le K_2^{\frac{3}{2}}$ in  $B_{\tau}(o,2r)\times[0,T]$, $\tau=T-t$. Let $u$ be a solution to the conjugate heat equation $\left(\Delta_t+\frac{\partial}{\partial t}-R \right)u=0$ in $B_{\tau}(o,2r)\times[0,T]$ satisfying $0<u\le A$. Then there exist nonnegative constants $\tilde{C_1}(n)$ and $\tilde{C_2}(n)$ such that
		\begin{equation}
			\begin{aligned}
				&\sup\limits_{x\in B_{\tau}(o,r)} \frac{|\nabla u|^2}{u^2}(x,\tau) \\
				&\le \left(1+\ln \frac{A}{u}\right)^2\left(\frac{1}{\tau}+\tilde{C_2}(n)K_1+2K_2+\frac{1}{r^2}(\tilde{C_1}(n)\sqrt{K_1}r\coth(\sqrt{K_1}r)+\tilde{C_2}(n))\right).\nonumber
			\end{aligned}
		\end{equation}
	\end{lemma}
	
	\begin{remark}\label{rm5.2} If the conditions $\Ric(x,t)\ge-K_1g(x,t)$, $|\nabla R|(x,t)\le K_2^{\frac{3}{2}}$ and $0<u\le A$ are satisfied on $M\times[0,T]$, then by letting $R\rightarrow \infty$, we have that
		\begin{equation}
			\frac{|\nabla u|^2}{u^2}(x,\tau) \le \left(1+\ln \frac{A}{u}\right)^2\left(\frac{1}{\tau}+\tilde{C_3}(n)(K_1+K_2)\right) \nonumber
		\end{equation}
		holds on $M\times(0,T]$.
	\end{remark}
	In \cite{KZ}, Kuang-Zhang established a Li-Yau type estimate for positive solutions of the conjugate heat equation under the Ricci flow.
	\begin{lemma}\label{lm5.3}(\cite{KZ})
		Let $M^n$ be an n-dimensional closed Riemannian manifold and $g(x,t), t\in [0,T)$ a solution to the Ricci flow  $\frac{\partial}{\partial t}g(x,t)=-2\Ric(x,t)$ on $M$. Assume that $R(x,t)\ge0$, $(x,t)\in M \times[0,T)$. Let $u$ be a positive solution to the conjugate heat equation $\left(\Delta_t+\frac{\partial}{\partial t}-R\right)u=0$ on  $ M\times [0,T)$. Let $u=(4\pi \tau)^{-\frac{n}{2}}e^{-f}$ with $\tau=T-t$. Then 
		\begin{equation}
			2\Delta f-|\nabla f|^2+R\le\frac{2n}{\tau},\nonumber
		\end{equation}
		or equivalently,
		\begin{equation}
			\frac{|\nabla u|^2}{u^2}+2\frac{u_{t}}{u}-R\le\frac{2n}{T-t}\nonumber
		\end{equation}
		holds on $ M \times[0,T)$.\\
	\end{lemma}

	\noindent{\it Proof of Theorem \ref{th1.9}.} Since the heat kernel $G(y,l;x,t)$ of the conjugate heat equation satisfies $G(y,l;x,t)=H(x,t;y,l)$, where $H(x,t; y,l)$ is the heat kernel of the heat equation and $x,y\in M$ and $0\le l<t \le T$, by Lemma \ref{lm3.2}, we have
	\begin{equation}\label{e5.1}
		\frac{\tilde{C_7}(n)e^{-\tilde{C_8}(n)e^{\tilde{C_9}(n)K_1T}}e^{-\frac{4d_T^2(M)}{T-l}}}{\vol_T(B_T(y,\sqrt{\frac{T-l}{8}}))}	\le G(y,l;x,T)\le \frac{\tilde{C_4}(n)e^{\tilde{C_5}(n)e^{\tilde{C_6}(n)K_1T}}}{\vol_l(B_l(x,\sqrt{\frac{T-l}{8}}))}. 
	\end{equation}
	For a time $0<s\le\frac{T}{2}$, we denote $h(y,l)=G(y,l;x,T)$, $x,y\in M$ and $l\in [0,s]$.\\
	Since $h(y,l)$ is a positive solution of the conjugate heat equation $\left(\Delta_{l,y}+\frac{\partial}{\partial l}-R\right)h(y,l)=0$, by Remark \ref{rm5.2}, we derive
	\begin{equation}\label{e5.4}
		|\nabla_y \ln G(y,l;x,T)|^2\le \left(1+\ln\frac{A}{G(y,l;x,T)}\right)^2\left(\frac{1}{l}+\tilde{C}_{10}(n)(K_1+K_2)\right),\ (y,l)\in M\times (0,s],
	\end{equation}
	where $A=\sup\limits_{(y,l)\in M\times [0,s]} G(y,l;x,T)$.
	
	The upper bound in \eqref{e5.1} gives
	\begin{equation}
		A\le \frac{\tilde{C_4}(n)e^{\tilde{C}_{11}(n)e^{\tilde{C}_{12}(n)K_1T}}}{\vol_0(B_0(x,\sqrt{\frac{T}{16}}))}.\nonumber
	\end{equation} 
	Hence, combining with the lower bound in \eqref{e5.1}, we have 
	\begin{equation}\label{e5.2}
		\begin{aligned}
			\frac{A}{G(y,l;x,T)} &\le \tilde{C}_{13}(n)e^{\tilde{C}_{14}(n)e^{\tilde{C}_{15}(n)K_1T}}e^{\frac{8d^2_T(M)}{T}}\frac{\vol_0(B_0(y,\sqrt{\frac{T}{8}}))}{\vol_0(B_0(x,\sqrt{\frac{T}{16}}))} \\
			& \le \tilde{C}_{13}(n)e^{\tilde{C}_{14}(n)e^{\tilde{C}_{15}(n)K_1T}}e^{\frac{8d^2_T(M)}{T}}\frac{\vol_0(B_0(x,d_0(M)+\sqrt{\frac{T}{8}}))}{\vol_0(B_0(x,\sqrt{\frac{T}{16}}))}. 
		\end{aligned}
	\end{equation}
	According to the Bishop-Gromov volume comparison theorem, we compute 
	\begin{equation}\label{e5.3}
		\frac{\vol_0(B_0(x,d_0(M)+\sqrt{\frac{T}{8}}))}{\vol_0(B_0(x,\sqrt{\frac{T}{16}}))} \le \left(\frac{d_0(M)+\sqrt{\frac{T}{8}} }{\sqrt{\frac{T}{16}}}\right)^ne^{\tilde{C}_{16}(n)\sqrt{K_1}(d_0(M)+\sqrt{T})}.
	\end{equation}
	Combining \eqref{e5.2} and \eqref{e5.3}, one has 
	\begin{equation}
		\frac{A}{G(y,l;x,T)} \le \tilde{C}_{13}(n)e^{\tilde{C}_{14}(n)e^{\tilde{C}_{15}(n)K_1T}}e^{\left[\frac{8d^2_T(M)}{T}+\tilde{C}_{17}(n)(\frac{d_0(M)}{\sqrt{T}}+1+d_0(M)\sqrt{K_1}+\sqrt{K_1T}) \right]}.\nonumber
	\end{equation}
	Plugging this into \eqref{e5.4} results in
	\begin{equation}
		|\nabla_y \ln G(y,l;x,T)|^2\le \tilde{C}_{20}(n)e^{\tilde{C}_{21}(n)K_1T} \left(1+\frac{d^2_0(M)}{T}\right)^2 \left(\frac{1}{l}+K_1+K_2\right),\ (y,l)\in M\times(0,s].\nonumber
	\end{equation}
	We take $l=s$ and use the arbitrariness of $0<s\le \frac{T}{2}$ to conclude 
	\begin{equation}
		|\nabla_y \ln G(y,l;x,T)|^2\le \tilde{C}_{20}(n)e^{\tilde{C}_{21}(n)K_1T} \left(1+\frac{d^2_0(M)}{T}\right)^2 \left(\frac{1}{l}+K_1+K_2\right),\ (y,l)\in M\times(0,\frac{T}{2}].\nonumber
	\end{equation}
	Since a positive solution $u(x,t)$ to the conjugate heat equation can be written as $u(x,t)=\int_M G(x,t;y,T)u(y,T)d\mu_T(y) $ with $t\in [0,T)$, the above estimate holds for any $u(x,t)$, i.e., 
	\begin{equation}
		|\nabla \ln u(x,t)|^2 \le \tilde{C}_{20}(n)e^{\tilde{C}_{21}(n)K_1T} \left(1+\frac{d^2_0(M)}{T}\right)^2 \left(\frac{1}{t}+K_1+K_2\right),\ (x,t)\in M\times(0,\frac{T}{2}].\nonumber
	\end{equation}
	Combining Lemma \ref{lm5.3}, we achieve that 
	\begin{equation}
		\begin{aligned}
			\frac{|\nabla u|^2}{u^2}+\frac{u_t}{u}&=\frac{1}{2}\left(\frac{|\nabla u|^2}{u^2}+2\frac{u_t}{u}-R\right)+\frac{1}{2}|\nabla \ln u|^2+\frac{1}{2}R \\
			&\le \frac{n}{T-t}+\tilde{C}_{22}(n)e^{\tilde{C}_{21}(n)K_1T} \left(1+\frac{d^2_0(M)}{T}\right)^2 \left(\frac{1}{t}+K_1+K_2\right),\nonumber 
		\end{aligned}
	\end{equation}
where $(x,t)\in M\times(0,\frac{T}{2}].$ This completes the proof of Theorem \ref{th1.9}.   \qed

\section*{Acknowledgements}

We would like to thank Prof. Qi S. Zhang for inspiring discussions and invaluable suggestions. We are also grateful for Prof. Huai-Dong Cao and Prof. Bin Qian for very helpful comments. Research is partially supported by NSFC Grant No. 11971168, Shanghai Science and Technology Innovation Program Basic Research Project STCSM 20JC1412900, and Science and Technology Commission of Shanghai Municipality (STCSM) No. 22DZ2229014.


\begin{thebibliography}{99}
		
		
		\bibitem{AB} Aronson, Donald G.; B\'enilan, Philippe, {\it R\'egularit\'e des solutions de l'\'equation des milieux poreux dans $R^N$}, C. R. Acad. Sci. Paris S\'er. A-B 288 (1979), no. 2, A103-A105.
		\bibitem{BCP} Bailesteanu, Mihai; Cao, Xiaodong; Pulemotov, Artem, {\it Gradient estimates for the heat equation under the Ricci flow}, J. Funct. Anal. 258 (2010), no. 10, 3517–3542. 
		\bibitem{BBG} Bakry, Dominique; Bolley, Fran\c{c}ois; Gentil, Ivan, {\it The Li-Yau inequality and applications under a curvature-dimension condition}, Ann. Inst. Fourier (Grenoble) 67 (2017), no. 1, 397-421. 
		\bibitem{BL}Bakry, Dominique; Ledoux, Michel,\emph{ A logarithmic Sobolev form of the Li-Yau parabolic inequality}, Rev. Mat. Iberoamericana, Volume 22, Number 2 (2006), 683-702.
		\bibitem{BQ} Bakry, Dominique; Qian, Zhongmin, {\it Harnack inequalities on a manifold with positive or negative Ricci curvature}, Rev. Mat. Iberoam. 15 (1999), no. 1, 143-179.
		\bibitem{CHD} Cao, Huai-Dong, {\it On Harnack's inequalities for the K\"ahler-Ricci flow}. Invent. Math. 109 (1992), no. 2, 247-263.
	        \bibitem{CN} Cao, Huai-Dong; Ni, Lei, {\it Matrix Li-Yau-Hamilton estimate for the heat equation on K\"ahler manifolds}, Math. Ann. 331 (2005), 795-807.
		\bibitem{CFL} Cao, Xiaodong; Fayyazuddin Ljungberg, Benjamin; Liu, Bowei, {\it Differential Harnack estimates for a nonlinear heat equation}, J. Funct. Anal. 265 (2013), no. 10, 2312–2330.
		\bibitem{CH}Cao, Xiaodong; Hamilton, Richard S., \emph{Differential Harnack estimates for time-dependent heat equations with potentials}, Geom. Funct. Anal. 19 (2009) 989–1000.
		\bibitem{CxZz} Cao, Xiaodong; Zhang, Zhou, {\it Differential Harnack estimates for parabolic equations}, Complex and differential geometry, 87–98, Springer Proc. Math., 8, Springer, Heidelberg, 2011.
		\bibitem{BLN} Chow, Bennett; Lu, Peng; Ni, Lei,  \emph{Hamilton’s Ricci Flow}, Graduate Studies in Mathematics, vol. 77, American Mathematical Society/Science Press, Providence, RI/New York, 2006.
		\bibitem{GM} Garofalo, Nicola; Mondino, Andrea, \emph{Li–Yau and Harnack type inequalities in $RCD^*(K; N)$ metric measure spaces}, Nonlinear Anal. 95 (2014) 721–734.
		\bibitem{HM}Hallgren, Max, \emph{The Entropy of Ricci Flows with Type-I Scalar Curvature Bounds}, arXiv:2007.10376.
		\bibitem{Ha}Hamilton, Richard S., \emph{A matrix Harnack estimate for the heat equation}. Comm. Anal. Geom. 1 (1993), no.1, 113-126.
		\bibitem{Ha2}Hamilton, Richard S., \emph{The Ricci flow on surfaces}, in: Mathematics and General Relativity, SantaCruz, CA, 1986, in: Contemp. Math., vol. 71, Amer. Math. Soc., Providence, RI, 1988, pp. 237–262.
		\bibitem{Ha4}Hamilton, Richard S., \emph{The Harnack estimate for the Ricci flow}, J. Differential Geom. 37 (1) (1993)225–243
		\bibitem{KZ}Kuang, Shilong , Zhang, Qi S.,\emph{ A gradient estimate for all positive solutions of the conjugate heat equation under Ricci flow}, J. Funct. Anal. 255 (4) (2008) 1008–1023.
		\bibitem{Lij}Li, Jiayu, {\it Gradient estimates and Harnack inequalities for nonlinear parabolic and nonlinear elliptic equations on Riemannian manifolds}, J. Funct. Anal. 100 (1991), no. 2, 233–256. 
		\bibitem{LX} Li, Junfang; Xu, Xiangjin, \emph{Differential Harnack inequalities on Riemannian manifolds I: linear heat equation},  Adv. Math. 226 (5) (2011) 4456–4491.
		\bibitem{LY}Li, Peter; Yau, Shing-Tung, \emph{On the parabolic kernel of the Schrödinger operator}. Acta Math. 156(1986), no. 3-4, 153–201.
		\bibitem{Liy}Li, Yi, \emph{Li-Yau-Hamilton estimates and Bakry-Emery-Ricci curvature}. Nonlinear Anal. 113 (2015), 1–32.
		\bibitem{LZ} Li, Xiaolong; Qi S. Zhang, {\it Matrix Li-Yau-Hamilton estimates under Ricci Flow and parabolic frequency}, arXiv:2306.10143.
        %\bibitem{Liu}Liu, Shiping, \emph{Gradient estimates for solutions of the heat equation under Ricci flow}, Pacific J. Math. 243:1 (2009), 165–180.
		\bibitem{Per} Perelman, Grisha, {\it The entropy formula for the Ricci flow and its geometric applications}, arXiv:math/0211159.
		\bibitem{Qian} Qian, Bin, {\it Remarks on differential Harnack inequalities}, J. Math. Anal. Appl. 409 (2014), no. 1, 556-566.
		\bibitem{QZZ} Qian, Zhongmin; Zhang, Hui-Chun; Zhu, Xi-Ping, \emph{Sharp spectral gap and Li–Yau’s estimate on Alexandrov spaces}, Math. Z. 273 (3–4) (2013) 1175–1195.	 
		\bibitem{WFY} Wang, Fengyu,  \emph{Gradient and Harnack inequalities on noncompact manifolds with boundary}, Pacific J. Math. 245 (1) (2010) 185–200.
		\bibitem{Yau} Yau, Shing-Tung, \emph{On the Harnack inequalities of partial differential equations}, Comm. Anal. Geom. 2 (1994), no. 3, p. 431-450.		
		\bibitem{YZ}Yu, Chengjie; Zhao, Feifei, \emph{Li-Yau multiplier set and optimal Li-Yau gradient estimate on hyperbolic spaces}. Potential Anal. 56 (2022), no. 2, 191–211.
		\bibitem{ZhZx} Zhang, Hui-Chun; Zhu, Xi-Ping, \emph{Local Li-Yau’s estimates on $RCD^*
(K, N)$ metric measure spaces}, Calc. Var. Part. Diff. Equ. 55 (2016), no. 4, Paper No. 93, 30 pp.
		\bibitem{ZQ}Zhang, Qi S., \emph{A Sharp Li-Yau gradient bound on Compact Manifolds }, arXiv:2110.08933.
		\bibitem{ZQ2}Zhang, Qi S., \emph{Some gradient estimates for the heat equation on domains and for an equation by Perelman}, Int. Math. Res. Not. (2006), Art. ID 92314, 39 pp.
		\bibitem{ZZ1}Zhang, Qi S.; Zhu, Meng,  \emph{Li-Yau gradient bound for collapsing manifolds under integral curvature condition}. Proc. Amer. Math. Soc. 145 (2017), no. 7, 3117–3126.
		\bibitem{ZZ2}Zhang, Qi S.; Zhu, Meng,  \emph{Li-Yau gradient bounds on compact manifolds under nearly optimal curvature conditions}.  J. Funct. Anal. 275 (2018), no. 2, 478–515.
		\bibitem{ZM} Zhu, Meng, \emph{Davies type estimate and the heat kernel bound under the Ricci flow}. Trans. Amer. Math. Soc. 368 (2016), no. 3, 1663–1680.
	\end{thebibliography}
\end{document}